\documentclass[12pt]{amsart}
\usepackage{mathrsfs}
\usepackage{amsfonts}
\usepackage{amssymb}
\usepackage{amscd, amsmath, amsthm, amsxtra}
\usepackage[papersize={7.7in,10.3in},textwidth=15.2cm,textheight=21cm,centering]{geometry}
\usepackage{enumerate}

\usepackage{xcolor}
\definecolor{cite}{rgb}{0.00,0.60,1.00}
\definecolor{url}{rgb}{1.00,0.10,0.80}
\definecolor{link}{rgb}{0.00,0.00,1.00}
\usepackage[colorlinks,linkcolor=link,urlcolor=url,citecolor=cite,pagebackref,breaklinks]{hyperref}

\hypersetup{
pdfstartpage=1,
pdfstartview=FitH}

\DeclareFontFamily{U}{mathx}{\hyphenchar\font45}
\DeclareFontShape{U}{mathx}{m}{n}{
      <5> <6> <7> <8> <9> <10>
      <10.95> <12> <14.4> <17.28> <20.74> <24.88>
      mathx10
      }{}
\DeclareSymbolFont{mathx}{U}{mathx}{m}{n}
\DeclareMathAccent{\widecheck}{\mathalpha}{mathx}{"71}

 \usepackage{caption} 
\numberwithin{equation}{section}

\allowdisplaybreaks

\newtheorem{theorem}{Theorem}[section]
\newtheorem{lemma}{Lemma}[section]

\newtheorem{proposition}{Proposition}[section]

\makeatletter
\newcounter{roem}
\renewcommand{\theroem}{\Roman{roem}}

\newcommand{\c@org@eq}{}
\let\c@org@eq\c@equation
\newcommand{\org@theeq}{}
\let\org@theeq\theequation

\newcommand{\setroem}{
\let\c@equation\c@roem
 \let\theequation\theroem}

\newcommand{\setarab}{
\let\c@equation\c@org@eq
\let\theequation\org@theeq}
\makeatother

\newtheorem*{claim*}{Claim}

\theoremstyle{remark}
\newtheorem{remark}{\bf Remark}

\newcommand{\ud}{\mathrm{d}}

\newcommand{\F}{\mathbb{F}}
\newcommand{\Q}{\mathbb{Q}}
\newcommand{\R}{\mathbb{R}}
\newcommand{\Z}{\mathbb{Z}}

\newcommand{\cA}{\mathcal{A}}
\newcommand{\cB}{\mathcal{B}}
\newcommand{\cC}{\mathcal{C}}

\newcommand{\cR}{\mathcal{R}}

\newcommand{\fA}{\mathfrak{A}}

\newcommand{\fM}{\mathfrak{M}}

\newcommand{\fS}{\mathfrak{S}}

\usepackage{graphicx}
\usepackage{tikz}

\begin{document}

\vglue -2mm

\title
{A conjecture of S\' ark\" ozy on quadratic residues. II}
\author{Yong-Gao Chen}
\address{School of Mathematical Sciences and Institute of Mathematics, Nanjing Normal University,  Nanjing  210023,  P. R. China}
\email{ygchen@njnu.edu.cn}

\author{Ping Xi }
\address{School of Mathematics and Statistics, Xi'an Jiaotong University, Xi'an 710049, P. R. China}
\email{ping.xi@xjtu.edu.cn}

\subjclass[2020]{11P70, 11B99}

\keywords{sumset, additive decomposition, 
quadratic residue}

\begin{abstract} 
Denote by $\cR_p$ the set of all quadratic residues in $\F_p$ for
each prime $p$. A conjecture of A. S\'ark\"ozy asserts, for all
sufficiently large $p$, that no subsets $\cA,\cB\subseteq\F_p$
with $|\cA|,|\cB|\geqslant2$ satisfy $\cA+\cB=\cR_p$. In this
paper, we show that if such subsets $\cA,\cB$ do exist, then there
are at least $(\log 2)^{-1}\sqrt p-1.6$ elements in $\cA+\cB$ that
have unique representations and one should have
\begin{align*}
\frac{1}{4}\sqrt{p}< |\cA|,|\cB|< 2\sqrt{p}-1.
\end{align*}
This refines previous bounds obtained by I.E. Shparlinski, I.D.
Shkredov, and Y.-G. Chen and X.-H. Yan. Moreover, we also establish bounds for $|\cA|,|\cB|$ and the
additive energy $E(\cA,\cB)$ if few elements in $\cA+\cB$ have
unique representations.
\end{abstract}

\maketitle

\section{Introduction}\label{sec:Introduction}

For each prime $p$, denote by $\F_p$ the finite field of $p$
elements and by $\cR_p$ the set of all quadratic residues modulo
$p$. For any subsets $\cA,\cB\subseteq\F_p,$ define the sumset
$\cA+\cB =\{a+b: a\in\cA, b\in \cB\}.$ A. S\'ark\"ozy \cite{Sa12}
conjectured that for all sufficiently large primes $p$, $\cR_p$
has no $2$-decomposition $\cA+\cB=\cR_p$ with
$|\cA|,|\cB|\geqslant2,$ which is also believed to be valid for
all odd primes $p$.

There are some partial results towards the above conjecture,
although the full generality seems beyond the current approach.
S\'ark\"ozy \cite{Sa12} considered the ternary analogue, and
proved that  for all sufficiently large primes $p$, no subsets
$\cA,\cB,\cC\subseteq\F_p$ with $|\cA|,|\cB|,|\cC|\geqslant2$
satisfy $\cA+\cB+\cC=\cR_p$. Recently, Y.-G. Chen and X.-H. Yan
\cite{CY21} proved this claim for all primes $p$. In addition, if
there do exist certain $\cA,\cB\subseteq\F_p$ such that
$|\cA|,|\cB|\geqslant2$ and $\cA+\cB=\cR_p,$ S\'ark\"ozy
\cite{Sa12} showed that
\begin{align}\label{eq:upper-lowerbound-Sarkozy}
\frac{\sqrt{p}}{3\log p}\leqslant |\cA|,|\cB|\leqslant
\sqrt{p}\log p.
\end{align}
The factor $\log p$ was shortly removed by I.E. Shparlinski
\cite{Shp13} and I.D. Shkredov \cite{Shk14} with some refined
constants. More remarkably, Shkredov \cite{Shk14} proved the above
conjecture of S\'ark\"ozy in the case $\cA=\cB,$ i.e.,
$\cA+\cA\neq\cR_p$ holds for any $\cA\subseteq\F_p$ with
$|\cA|\geqslant2.$ A quantitative refinement on
\eqref{eq:upper-lowerbound-Sarkozy} was recently given by Y.-G.
Chen and X.-H. Yan \cite{CY21}, who proved that
\begin{align}\label{eq:upper-lowerbound-ChenYan}
\frac{7-\sqrt{17}}{16}\sqrt{p}+1\leqslant |\cA|,|\cB|\leqslant
\frac{7+\sqrt{17}}{4}\sqrt{p}-6.63.
\end{align}
Note that $(7-\sqrt{17})/16\approx0.1798$ and
$(7+\sqrt{17})/4\approx2.7808$. This improves a previous result by
Shkredov \cite{Shk14}, in which the constants in the lower and
upper bounds are $1/6-o(1)$ and $3+o(1)$ as $p\rightarrow+\infty.$
The main tool of Shparlinski \cite{Shp13}  includes an estimate of
the double character sum
\[\sum_a\sum_b\chi(a+b)\]
due to A.A. Karatsuba \cite{Ka92}, where $\chi$ is a non-trivial
multiplicative character in $\F_p^\times.$ Alternatively, Shkredov
\cite{Shk14} and Chen and Yan \cite{CY21} employed the Weil bound
for complete character sums over finite fields; see Lemmas
\ref{lm:Weil} and \ref{lm:Wan} below.

In this paper, we make further refinements on the constants in
\eqref{eq:upper-lowerbound-ChenYan}.

\begin{theorem}\label{thm:AB-bound-explicit}
Let $p$ be an odd prime. Assume that $\cA+\cB=\cR_p$ is a
$2$-decomposition with $|\cA|,|\cB|\geqslant 2.$ Then
\begin{align}\label{eq:AB-bound-explicit}
\frac{1}{4}\sqrt{p}+\frac 18\leqslant |\cA|,|\cB|<
2\sqrt{p}-1.
\end{align}
\end{theorem}

Our next results involve the representation function
$r=r_{\cA+\cB}$ defined by
\begin{align*}
r(x)=r_{\cA+\cB}(x):=|\{(a,b)\in\cA\times\cB:a+b=x\}|.
\end{align*}
For $\theta\in\R,$ define the moment of $r$ by
\begin{align*}
M_\theta:=\sum_{x\in\cA+\cB}r(x)^\theta.
\end{align*}
Note that $M_0=|\cA+\cB|$, $M_1=|\cA||\cB|$. Moreover,
$M_2=E(\cA,\cB),$ the additive energy of $\cA,\cB,$ i.e.,
\begin{align*}
E(\cA,\cB)=|\{(a_1,a_2,b_1,b_2)\in\cA\times\cA\times\cB\times\cB:a_1+b_1=a_2+b_2\}|.
\end{align*}
We also introduce
\begin{align*}
\fA(\cA,\cB):=|\{x\in\cA+\cB:r(x)=1\}|
\end{align*}
to characterize how many elements in $\F_p$ have unique
representations in $\cA+\cB.$ Clearly,
\begin{align*}
0\leqslant \fA(\cA,\cB)\leqslant |\cA+\cB|\leqslant |\cA||\cB|\leqslant E(\cA,\cB).
\end{align*}
Roughly speaking, the three quantities $|\cA+\cB|$, $|\cA||\cB|$
and $E(\cA,\cB)$ should be of similar sizes if $\fA(\cA,\cB)$ is
reasonably large; on the other hand,  $|\cA||\cB|$ becomes closer
to $2|\cA+\cB|$ as $\fA(\cA,\cB)/|\cA+\cB|$ decreases to zero (see
Lemma \ref{lm:AB,A+B-bound} for details). It is clear that if
$\fA(\cA,\cB) =|\cA+\cB|$, then
\begin{align*}
 \fA(\cA,\cB)= |\cA+\cB|= |\cA||\cB|=E(\cA,\cB).
\end{align*}

We are now ready to state our results.

\begin{theorem}\label{thm:uniquerepresentation}
Let $p$ be a prime. Assume that $\cA+\cB=\cR_p$ is a
$2$-decomposition with $|\cA|,|\cB|\geqslant 2.$ Then
\begin{align*}
\fA(\cA,\cB)\geqslant (\log 2)^{-1} \sqrt{p}-1.6.
\end{align*}
\end{theorem}

\begin{theorem}\label{thm:energy,uniquerepresentation}
Let $p$ be a prime. Assume that $\cA+\cB=\cR_p$ is a
$2$-decomposition with $|\cA|,|\cB|\geqslant 2.$ If there exists
some $\eta\in[0,\frac{1}{2})$ such that
\begin{align*}
\fA(\cA,\cB)\leqslant \eta (p-1),
\end{align*}
then we have
\begin{align*}
E(\cA,\cB)\geqslant\Big(\frac{1}{2}+\frac{2^{2-4\eta}-2^{3-2\eta}\eta+4\eta-1}{2-4\eta}\Big)(p-1)
\end{align*}
and
\begin{align*}
|\cA|,|\cB|\geqslant \frac{1}{2\cdot 4^\eta}\sqrt{p}.
\end{align*}
\end{theorem}

Theorem \ref{thm:uniquerepresentation} illustrates that there are suitably many elements in $\F_p$ which can be
expressed as $a+b$ with $(a,b)\in\cA\times\cB$ in a unique way, and however, if there are not too many such elements,  Theorem \ref{thm:energy,uniquerepresentation} shows that $\cA,\cB$
should be of reasonably large sizes.

In principle, there are two novelties in this paper. On one hand, the improvements revealed by Theorem \ref{thm:AB-bound-explicit} rely on Lemma \ref{lm:Wan} on estimates for complete character sums, which allows us to refine the constant derived directly from the Weil bound in the usual shape. In the last section, we will also outline the robustness of this bound in view of equidistributions of Birch \cite{Bi68} and Katz and Sarnak \cite{KS99}.
On the other hand, a new quantity $\fA(\cA,\cB)$ has been introduced to characterize the number of elements in $\F_p$ which have unique representations in $\cA+\cB$. We hope this should be of independent interests and will receive more attentions in further researches.

By closing this section, we remark that the arguments in this paper admit the bounds
\begin{align*}
E(\cA,\cB)\leqslant p^{3/2}-p,
\end{align*}
while the trivial bound would be $E(\cA,\cB)\leqslant |\cA||\cB|\min\{|\cA|,|\cB|\}$. We do not know how to beat the exponent $3/2.$

\smallskip

\section{Auxiliary results}\label{sec:lemmas}
We start with the following estimate for complete character sums
over finite fields due to A. Weil \cite{We48}.

\begin{lemma}\label{lm:Weil}
Let $\chi$ be a multiplicative character of order $d>1$ of $\F_p$.
Assume that $f(X)\in\F_p[X]$ has $k$ distinct zeros in the
algebraic closure of $\F_p$ and it is not a constant multiple of
the $d$-th power of a polynomial over $\F_p$. Then
\begin{align*}
\Big|\sum_{x\in\F_p}\chi(f(x))\Big|\leqslant (k-1)\sqrt{p}.
\end{align*}
\end{lemma}

Denote by $\chi_2$ the quadratic multiplicative character of
$\F_p^\times$, i.e., the Legendre symbol mod $p$. For
$\mathbf{a}=(a_1, \cdots , a_k)\in\F_p^k$, define
\begin{align*}
\fS_k(\mathbf{a};p)=\sum_{x\in\F_p} \chi_2((x+a_1)\cdots (x+a_k)).
\end{align*}

For $k=2$, an elementary argument leads to
\begin{align*}
\fS_2(\mathbf{a};p)=-1
\end{align*}
if $a_1\neq a_2,$ and $\fS_2(\mathbf{a};p)=p-1$ otherwise. For $k\geqslant3$, it follows from Lemma \ref{lm:Weil} that
\begin{align}\label{eq:Weil}
|\fS_k(\mathbf{a};p)|\leqslant(k-1)\sqrt{p}
\end{align}
if the coordinates $a_1,\cdots, a_k$ are pairwise distinct. As in Chen and Yan \cite{CY21}, define
$c_k(p)$ to be the smallest real number $c$ such that
\begin{align*}
\fS_k(\mathbf{a};p)\leqslant c\sqrt{p}
\end{align*}
holds for all pairwise distinct elements $a_1, \cdots , a_k\in\F_p
$. Let $c_k$ be the supremum  of $c_k(p)$ over odd primes $p$.
The following result gives an upper bound for $c_k$, which is in fact a special case of D. Wan \cite[Corollary 2.3]{Wa97}.

\begin{lemma}[Wan \cite{Wa97}]\label{lm:Wan}
Let $k>0$ be an even integer. With the above notation, we have
\begin{align}\label{eq:ck-upperbound}
c_k\leqslant k-2.
\end{align}
More precisely, for each odd prime $p$ and all pairwise distinct
elements $a_1, \cdots , a_{k}\in\F_p,$ we have
\begin{align}\label{eq:Wan}
\fS_k(\mathbf{a};p)\leqslant(k-2)\sqrt{p}-1.
\end{align}
\end{lemma}

\proof To convince the readers, we now give a completely
elementary argument, plus the original Weil bound \eqref{eq:Weil},
that leads to \eqref{eq:Wan}.

Making the shift $x\rightarrow x-a_1,$ we may write
\begin{align*}
\fS_k(\mathbf{a};p) &=\sum_{x\in\F_p^\times}
\chi_2(x(x+a_2-a_1)\cdots (x+a_{k}-a_1)).
\end{align*}
Note that $\chi_2(x+h)=\chi_2(x)\chi_2(1+h\overline{x})$ for all
$x\in\F_p^\times.$ We then derive that
\begin{align*}
\fS_k(\mathbf{a};p)
&=\sum_{x\in\F_p^\times} \chi_2((1+(a_2-a_1)\overline{x})\cdots ((1+(a_{k}-a_1)\overline{x}))\\
&=\sum_{x\in\F_p} \chi_2((1+(a_2-a_1)x)\cdots
((1+(a_{k}-a_1)x))-1.
\end{align*}
Now the desired inequality follows from Lemma \ref{lm:Weil}
immediately.
\endproof

\begin{remark}
Since $c_k(p)$ is utilized to bound $\fS(\mathbf{a};p)$ for all pairwise distinct
elements $a_1, \cdots , a_{k}\in\F_p,$ it is not possible to improve \eqref{eq:ck-upperbound} with a smaller constant. This will be discussed in detail in the last section.
\end{remark}

The following lemma shows an auxiliary bound involving cardinalities
of $\cA$ and $\cB,$ which will be our starting point to prove
Theorem \ref{thm:AB-bound-explicit}.

\begin{lemma}\label{lm:AB-bound}
Assume that $\cA,\cB$ are two subsets of $\F_p$ with $
\cA+\cB\subseteq \cR_p,$  then we have
\begin{align}\label{eq:AB-bound}
|\cB||\cA| (|\cA|+2)^2\leqslant & 2\sqrt{p}(|\cA|^2-1)(|\cA|-2)\nonumber\\
& -|\cA|^3+11|\cA|-15+p(3|\cA|+2).
\end{align}
\end{lemma}

\proof The proof goes quite similarly to \cite[Lemma 2.7]{CY21}.
Consider the moment
\begin{align*}
\fM:=\sum_{b\in\cB}\Big(\sum_{a\in\cA}\chi_2(a+b)\Big)^2
\Big(\sum_{a\in\cA}\chi_2(a+b) +2 \Big)^2 .
\end{align*}
Since $\cA+\cB\subseteq \cR_p,$ we find that $\chi_2(a+b)=1$ for
all $a\in\cA,b\in\cB,$ which gives
\begin{align}\label{eq:M=A^4B}
\fM=|\cB||\cA|^2(|\cA|+2)^2.
\end{align}
On the other hand,
\begin{align*}
\fM&\leqslant
\sum_{x\in\F_p}\Big(\sum_{a\in\cA}\chi_2(x+a)\Big)^2\Big(\sum_{a\in\cA}\chi_2(a+b)
+2 \Big)^2\\
&=\sum_{\mathbf{a}\in\cA^4}\fS_4(\mathbf{a};p)
+4\sum_{\mathbf{a}\in\cA^3}\fS_3(\mathbf{a};p)+4\sum_{\mathbf{a}\in\cA^2}\fS_2(\mathbf{a};p).
\end{align*}
We may classify all $\mathbf{a}\in\cA^j$ $(j=2,3,4)$ according to
the multiples in coordinates. From  Lemma \ref{lm:Weil} and Lemma
\ref{lm:Wan}, it follows that
\begin{align*}
\fM \leqslant & (2\sqrt{p} -1)|\cA|(|\cA|-1)(|\cA|-2)(|\cA|-3)-6|\cA|(|\cA|-1)(|\cA|-3)\\
& -4|\cA|(|\cA|-1)+3(p-2)|\cA|(|\cA|-1)+(p-1)|\cA|\\
& +8\sqrt{p} |\cA|(|\cA|-1)(|\cA|-2) +12|\cA|(|\cA|-1)\\
&-4|\cA|(|\cA|-1)+4(p-1) |\cA|,
\end{align*}
from which  and \eqref{eq:M=A^4B} the lemma follows.
\endproof

\begin{lemma}\label{lm:p>AB}
Assume that $\cA,\cB$ are two subsets of $\F_p$ with $
\cA+\cB\subseteq \cR_p,$ then
\begin{align*}
p|\cA||\cB|\leqslant & (p-|\cA|) (p-|\cB|).
\end{align*}
\end{lemma}

\proof  Using the orthogonality of additive characters, by $
\cA+\cB\subseteq \cR_p,$ we derive that
\begin{align*}
|\cA||\cB|
&=\sum_{x\in\F_p}\chi_2(x)\mathop{\sum_{a\in\cA}\sum_{b\in\cB}}_{a+b=x}1\\
&=\frac{1}{p}\sum_{\psi\in\widehat{\F}_p}\Big(\sum_{a\in\cA}\psi(a)\Big)\Big(\sum_{b\in\cB}\psi(b)\Big)\Big(\sum_{x\in\F_p}\chi_2(x)\psi(-x)\Big),
\end{align*}
where $\widehat{\F}_p$ denotes the group of additive characters of
$\F_p.$ The last sum over $x$ vanishes if $\psi$ is trivial and is
of modulus $\sqrt{p}$ otherwise. Therefore,
\begin{align*}
|\cA||\cB|
&\leqslant\frac{1}{\sqrt{p}}\sum_{1\neq\psi\in\widehat{\F}_p}\Big|\sum_{a\in\cA}\psi(a)\Big|\Big|\sum_{b\in\cB}\psi(b)\Big|,
\end{align*}
from which and Cauchy's inequality we find
\begin{align*}
p(|\cA||\cB|)^2
&\leqslant\sum_{1\neq\psi\in\widehat{\F}_p}\Big|\sum_{a\in\cA}\psi(a)\Big|^2\times
\sum_{1\neq\psi\in\widehat{\F}_p}\Big|\sum_{b\in\cB}\psi(b)\Big|^2.
\end{align*}
Note that
\begin{align*}
\sum_{1\neq\psi\in\widehat{\F}_p}\Big|\sum_{a\in\cA}\psi(a)\Big|^2=\sum_{\psi\in\widehat{\F}_p}\Big|\sum_{a\in\cA}\psi(a)\Big|^2-|\cA|^2=p|\cA|-|\cA|^2,
\end{align*}
and a similar identity can also be derived for the other sum. We
now arrive at
\begin{align*}
p(|\cA||\cB|)^2 &\leqslant(p|\cA|-|\cA|^2)(p|\cB|-|\cB|^2),
\end{align*}
from which Lemma \ref{lm:p>AB} follows.
\endproof

The following lemma is taken directly from \cite[Lemma 2.6]{CY21}.
\begin{lemma}\label{lm:A>4}
Assume that there exist subsets $\cA,\cB\subseteq\F_p$ with
$|\cB|\geqslant |\cA|\geqslant 2$ and $ \cA+\cB=\cR_p,$ then
$|\cA|\geqslant 5$.
\end{lemma}

\smallskip

\section{Proof of Theorem \ref{thm:AB-bound-explicit}}

Firstly we deal with the lower bound in
\eqref{eq:AB-bound-explicit}.
We now assume $\cA,\cB$ form a $2$-decomposition of  $\cR_p,$ so
that
\begin{align}\label{eq:AB>=p/2}
|\cA||\cB|\geqslant |\cA+\cB|=\frac{p-1}{2}.
\end{align}
By Lemma \ref{lm:A>4}, $|\cA|\geqslant 5$. On the other hand, it follows
from Lemma \ref{lm:AB-bound} that
\begin{align*}
|\cB||\cA| (|\cA|+2)^2\leqslant &
2\sqrt{p}(|\cA|^2-1)(|\cA|-2)\\
&-|\cA|^3+11|\cA|-15+p(3|\cA|+2).
\end{align*}
Therefore,
\begin{align*}
\frac{p-1}{2}(|\cA|+2)^2\leqslant &
2\sqrt{p}(|\cA|^2-1)(|\cA|-2)\\
&-|\cA|^3+11|\cA|-15+p(3|\cA|+2),
\end{align*}
from which and $|\cA|\geqslant 5$  we obtain
\begin{align*}
p(|\cA|^2-2|\cA|)\leqslant & 4\sqrt{p}(|\cA|^2-1)(|\cA|-2)
-2|\cA|^3+|\cA |^2+26|\cA|-26
\\
\leqslant & 4\sqrt{p}|\cA|(|\cA|^2-2|\cA|) .
\end{align*}
It follows that $ \sqrt{p}\leqslant 4|\cA|$. Now the lower bound
in \eqref{eq:AB-bound-explicit} follows from
\begin{align*}
\frac{1}{4}\sqrt{p} \leqslant
&\frac{(|\cA|^2-1)(|\cA|-2)}{|\cA|^2-2|\cA|}
-\frac{2|\cA|^3-|\cA |^2-26|\cA|+26}{4\sqrt{p}(|\cA|^2-2|\cA|) } \\
\leqslant &|\cA|-\frac{1}{|\cA|} -\frac{2|\cA|^3-|\cA |^2-26|\cA|+26}{16 |\cA|(|\cA|^2-2|\cA|) } \\
< & |\cA|-\frac 18.
\end{align*}

We now turn to prove the upper bound in
\eqref{eq:AB-bound-explicit} by contradiction. Suppose that $ |\cB|
\geqslant 2\sqrt{p}-1$. Then by Lemma \ref{lm:p>AB} we have
$$p\geqslant |\cA| (|\cB|+1) \geqslant 2\sqrt{p} |\cA|.$$
It follows that $\sqrt{p} \geqslant 2 |\cA|$.  The lower bound in
\eqref{eq:AB-bound-explicit} yields $\sqrt{p}<4|\cA|-\frac 12$.
From Lemma \ref{lm:AB-bound} it follows that
\begin{align*}
|\cB||\cA| (|\cA|+2)^2\leqslant &
2\sqrt{p}(|\cA|^2-1)(|\cA|-2)\\
&-|\cA|^3+11|\cA|-15+p(3|\cA|+2)\\
=& 2\sqrt{p}|\cA| (|\cA|+2)^2 \\
&+p(3|\cA|+2)- 2\sqrt{p} (6 |\cA|^2+5|\cA|-2)-|\cA|^3+11|\cA|-15.
\end{align*}
Since $2 |\cA|\leqslant \sqrt{p} <4|\cA|-\frac 12$, it follows
that
\begin{align*}
&p(3|\cA|+2)- 2\sqrt{p} (6 |\cA|^2+5|\cA|-2)\\
\leqslant & \max\Big\{ (2 |\cA|)^2 (3|\cA|+2)- 2(2 |\cA|) (6
|\cA|^2+5|\cA|-2),\\
& \left(4|\cA|-\frac 12\right)^2 (3|\cA|+2)- 2\left(4|\cA|-\frac 12\right) (6 |\cA|^2+5|\cA|-2)\Big\} \\
=&\left(4|\cA|-\frac 12\right)^2(3|\cA|+2)- 2\left(4|\cA|-\frac
12\right) (6
|\cA|^2+5|\cA|-2)\\
=& -14 |\cA|^2 +\frac{55}4 |\cA|-\frac32.
\end{align*}
Hence
\begin{align*}
|\cB||\cA| (|\cA|+2)^2\leqslant & 2\sqrt{p}|\cA| (|\cA|+2)^2 -14
|\cA|^2 +\frac{55}4
|\cA|-\frac32-|\cA|^3+11|\cA|-15\\
=&2\sqrt{p}|\cA| (|\cA|+2)^2 -|\cA|^3-14 |\cA|^2 +\frac{99}4
|\cA|-\frac{33}2 .
\end{align*}
In view of $|\cA|\geqslant 5$,
$$-|\cA|^3-14 |\cA|^2 +\frac{99}4 |\cA|-\frac{33}2 < -|\cA|
(|\cA|+2)^2.$$ So
$$|\cB||\cA| (|\cA|+2)^2<2\sqrt{p}|\cA| (|\cA|+2)^2-|\cA|
(|\cA|+2)^2.$$ This gives a contradiction to $ |\cB| \geqslant
2\sqrt{p}-1$, and thus completes the proof of Theorem
\ref{thm:AB-bound-explicit}.

\smallskip

\section{Proofs of Theorems \ref{thm:uniquerepresentation} and \ref{thm:energy,uniquerepresentation}}
We start with the following lemma.
\begin{lemma}\label{lm:AB,A+B-bound}
For any $\theta>0,$ we have
\begin{align}\label{eq:AB,kappa-1}
2^\theta|\cA+\cB|^{\theta+1}\leqslant
(|\cA||\cB|)^\theta(|\cA+\cB|+(2^\theta-1)\fA),
\end{align}
and
\begin{align}\label{eq:AB,kappa-2}
|\cA||\cB|\leqslant \fA+\sqrt{(E(\cA,\cB)-\fA)(|\cA+\cB|-\fA)},
\end{align}
where we write $\fA=\fA(\cA,\cB).$
\end{lemma}

\begin{remark}
The first inequality illustrates that $|\cA||\cB|$ becomes closer
to $2|\cA+\cB|$ as $\fA /|\cA+\cB|$ decreases to zero. Noting that
$\fA \leqslant |\cA+\cB| \leqslant |\cA||\cB|$, by
\eqref{eq:AB,kappa-2} we have
$$ |\cA+\cB| \leqslant |\cA||\cB|\leqslant |\cA+\cB|
+\sqrt{(E(\cA,\cB)-\fA)(|\cA+\cB|-\fA)}.$$ It follows that
  $|\cA||\cB|$ approaches to $|\cA+\cB|$ as $\fA$
grows, as long as the additive energy $E(\cA,\cB)$ is under
control.
\end{remark}

\proof  From H\" older's inequality, for any $\theta>0$ we have
\begin{align*}
M_0\leqslant M_{1}^{\frac {\theta }{\theta +1}} M_{-\theta}^{\frac
1{\theta +1}}.
\end{align*}
It follows that
\begin{align}\label{eq:M0M1M-1}
2^\theta  M_0^{\theta +1}\leqslant M_{1}^{\theta } 2^\theta
M_{-\theta} .
\end{align}
Moreover,  we also find
\begin{align*}
M_{-\theta}-\fA&=\sum_{\substack{x\in\cA+\cB\\ r(x)\geqslant2}}
r(x)^{-\theta} \leqslant 2^{-\theta}\sum_{\substack{x\in\cA+\cB\\
r(x)\geqslant2}}1 =\frac{M_0-\fA}{2^\theta}.
\end{align*}
Thus,
\begin{align*}
2^\theta M_{-\theta}  \leqslant M_0+(2^\theta -1) \fA ,
\end{align*}
from which and \eqref{eq:M0M1M-1}, we may conclude
\eqref{eq:AB,kappa-1}.

Cauchy's inequality yields
\begin{align*}
(M_1-\fA)^2=\Big(\sum_{\substack{x\in\cA+\cB\\
r(x)\geqslant2}}r(x)\Big)^2 \leqslant
\Big(\sum_{\substack{x\in\cA+\cB\\
r(x)\geqslant2}}r(x)^2\Big)\Big(\sum_{\substack{x\in\cA+\cB\\
r(x)\geqslant2}}1\Big),\end{align*} from which it follows that
\begin{align*}
(|\cA||\cB|-\fA)^2\leqslant (E(\cA,\cB)-\fA)(|\cA+\cB|-\fA).
\end{align*}
This proves \eqref{eq:AB,kappa-2}, and thus completes the proof of
Lemma \ref{lm:AB,A+B-bound}.
\endproof

\begin{lemma}\label{lem:ABbound} If
$\fA(\cA,\cB)= \tau |\cA+\cB|$, then
\begin{align*}
|\cA||\cB|\geqslant 2^{1-\tau } |\cA+\cB|.
\end{align*}
\end{lemma}

\proof From \eqref{eq:AB,kappa-1} it follows that
\begin{align*}
|\cA||\cB|\geqslant
\frac{2|\cA+\cB|^{1+\frac{1}{\theta}}}{(|\cA+\cB|+(2^\theta-1)\tau
|\cA+\cB|)^{\frac{1}{\theta}}}=2(1+(2^\theta-1)\tau
)^{-\frac{1}{\theta}}\cdot |\cA+\cB|.
\end{align*}
Taking $\theta\rightarrow0^+,$ we arrive at
\begin{align*}
|\cA||\cB|\geqslant 2^{1-\tau } |\cA+\cB|.
\end{align*}
This completes the proof of Lemma \ref{lem:ABbound}.
\endproof

\proof[Proof of Theorem \ref{thm:uniquerepresentation}] Suppose
that $\fA(\cA,\cB)=\kappa (p-1)$. Since $\cA+\cB=\cR_p$ with
$|\cA|, |\cB| \geqslant 2$, it follows that $p\geqslant 37$ in
view of Lemmas \ref{lm:p>AB} and \ref{lm:A>4}. From Lemma
\ref{lm:p>AB}, we obtain
\begin{align*}
(p-1)|\cA||\cB| \leqslant p(p-|\cA|-|\cB|) \leqslant
p(p-2\sqrt{|\cA||\cB|}).
\end{align*}
This gives
\begin{align*}
|\cA||\cB|\leqslant\Big(\frac{p}{\sqrt{p}+1}\Big)^2.
\end{align*}
On the other hand, Lemma \ref{lem:ABbound} shows
\begin{align}\label{eq:ABlowerbound-kappa}
|\cA||\cB|\geqslant 2^{1-2\kappa} |\cA+\cB|=
\frac{1}{4^{\kappa}}\cdot (p-1).
\end{align}
We hence infer
\begin{align*}
\frac 1{4^{\kappa}}\cdot (p-1)\leqslant \Big(\frac{p}{\sqrt{p}+1}\Big)^2,
\end{align*}
from which it follows that
\begin{align*}
\kappa\log2\geqslant \frac{1}{2}\log(p-1)+\log(\sqrt{p}+1)-\log p\geqslant \frac{\sqrt{p}}{p-1}-\frac{1.6\log2}{p-1}
\end{align*}
for $p\geqslant37.$
This completes the proof of Theorem \ref{thm:uniquerepresentation} in view of $\fA(\cA,\cB)=\kappa(p-1)$.
\endproof

\proof[Proof of Theorem \ref{thm:energy,uniquerepresentation}] We
also assume that $\cA+\cB=\cR_p$ with $|\cA|, |\cB| \geqslant 2$
and write $\fA=\kappa (p-1)$. By  Theorem
\ref{thm:AB-bound-explicit} and \eqref{eq:ABlowerbound-kappa},
 we obtain
\begin{align*}
|\cA|(2\sqrt p-1)\geqslant |\cA||\cB|\geqslant
\frac{1}{4^{\kappa}}\cdot (p-1),
\end{align*}
giving
\begin{align*}
|\cA|\geqslant \frac{1}{2\cdot 4^{\kappa}}\cdot \frac{p-1}{\sqrt
p-1/2} \geqslant \frac{1}{2\cdot 4^{\kappa}}\cdot \sqrt{p} .
\end{align*}
A similar lower bound for $|\cB|$ also holds by symmetry.

On the other hand, we infer from \eqref{eq:AB,kappa-2} that
\begin{align*}
|\cA||\cB|\leqslant (\kappa+\sqrt{(E(\cA,\cB)/(p-1)-\kappa
)(1/2-\kappa)})(p-1).
\end{align*}
This, together with \eqref{eq:ABlowerbound-kappa}, yields
\begin{align*}
E(\cA,\cB)\geqslant\Big(\frac{1}{2}+\frac{2^{2-4\kappa}-2^{3-2\kappa}\kappa+4\kappa-1}{2-4\kappa}\Big)(p-1)
\end{align*}
as desired.

\smallskip

\section{Concluding remarks}

\subsection{More bounds for possible decompositions}

Theorem \ref{thm:AB-bound-explicit} shows that $|\cA|/|\cB|$
should be between $\frac{1}{8}$ and $8$ if $\cA+\cB=\cR_p$ is a
$2$-decomposition with $|\cA|,|\cB|\geqslant 2.$
In fact, the above arguments may give better bounds if the size of $|\cA|/|\cB|$ is taken into account.

\begin{proposition}\label{prop}
Let $p$ be a large prime. Assume that $\cA+\cB=\cR_p$ is a
$2$-decomposition with $|\cB|\geqslant|\cA|\geqslant 2$ and
$|\cA|=\delta|\cB|$ with $\delta\in (\frac{1}{8},1].$ Then
\begin{align*}
\sqrt{\delta (p-1)/2} \leqslant |\cA|\leqslant \min\{ 2\delta,
\sqrt{\delta} \} \sqrt{p-1},
\end{align*}
and
\begin{align*}
\sqrt{(p-1)/(2\delta)}\leqslant |\cB|\leqslant \min\{ 2,
\sqrt{\delta^{-1}} \} \sqrt{p-1}.
\end{align*}
\end{proposition}

In fact, from \eqref{eq:AB>=p/2} it
follows that
\begin{align}\label{eq:AB-lowerbound1}
|\cA|\geqslant \sqrt{\delta (p-1)/2},\quad |\cB|\geqslant \sqrt{
(p-1)/(2\delta )}.
\end{align}
By Lemma \ref{lm:p>AB} we have $|\cB||\cA|\leqslant p-1$. Thus,
\begin{align}\label{eq:AB-upperbound2}
|\cA|\leqslant \sqrt{\delta (p-1)},\quad |\cB|\leqslant \sqrt{
(p-1)/\delta }.
\end{align}
On the other hand, we infer from Lemma \ref{lm:AB-bound} that
\begin{align*}
|\cB|&\leqslant 2\sqrt{p}-1\leqslant 2\sqrt{p-1}
\end{align*}
and also
\begin{align*}
|\cA|=\delta|\cB|&\leqslant 2\delta \sqrt{p-1}.
\end{align*}
Hence Proposition \ref{prop} follows by collecting the above two inequalities.

\smallskip

In particular, if $\cA+\cB=\cR_p$ is a $2$-decomposition with
$|\cA|=|\cB|\geqslant2,$ then
\begin{align*}
\sqrt{\frac{p-1}{2}} \leqslant |\cA|=|\cB|\leqslant\sqrt{p-1}
\end{align*}
as in Proposition \ref{prop}. The constant
$1/\sqrt{2}\approx 0.7071$ is much larger than $1/4=0.25$ as in Theorem \ref{thm:AB-bound-explicit}.

\subsection{Comments on the character sum $\fS(\mathbf{a};p)$}
Recall the quantities $c_k(p)$ and $c_k$ defined in Section \ref{sec:lemmas} to bound the complete character sum $\fS(\mathbf{a};p)$.
It is clear that $c_2(p)=-1/\sqrt{p}$, and the determination of $c_k(p)$ is no longer elementary for $k\geqslant4.$ Here is our observation.
\begin{proposition}
Let $k\geqslant2$ be a fixed even integer.
For all sufficiently large primes $p,$ we have
\[c_k(p)=k-2+o(1).\]
In particular, we have $c_k=k-2.$
\end{proposition}

The fact that $c_4(p)=2+o(1)$ follows from a ``vertical" version of Sato--Tate conjecture for elliptic curves over finite fields as proven by B.J. Birch \cite{Bi68}. In general, the value of $c_k(p)$ corresponds to a related equidistribution of hyperelliptic curves of genus $g:=k/2-1$
defined by
\begin{align*}
E:y^2=(1+(a_2-a_1)x)\cdots((1+(a_{k}-a_1)x)
\end{align*}
with pairwise distinct $a_1,a_2,\cdots,a_k\in\F_p.$ N.M. Katz and P. Sarnak \cite{KS99} proved an equidistribution for hyperelliptic curves (in a more general setting) in families. In fact, they showed that the corresponding Frobenius classes become equidistributed in the unitary symplectic group $\mathrm{USp}(k-2)$ as $p\rightarrow+\infty$. More precisely, as long as $(a_1,a_2,\cdots,a_k)\in\F_p^k$ runs over all possibilities with pairwise distinct coordinates, the {\it normalized} character sum
\begin{align*}
\frac{1}{\sqrt{p}}\sum_{x\in\F_p}\chi_2((1+(a_2-a_1)x)\cdots((1+(a_{k}-a_1)x))
\end{align*}
becomes equidistributed in $[2-k,k-2]$, as $p\rightarrow+\infty$, with respect to the Haar measure with image in $[0,2\pi)^g$ given by
\begin{align*}
\frac{1}{g!}\Big(\prod_{1\leqslant i<j\leqslant g}(2\cos\theta_i-2\cos\theta_j)^2\Big)\cdot\prod_{1\leqslant j\leqslant g}\frac{2}{\pi}\sin^2\theta_j\ud\theta_j.
\end{align*}

Keeping the uniformality in $\mathbf{a}\in\F_p^k,$ it is not possible to improve the constant $k-2$ in Lemma \ref{lm:Wan}.
On the other hand, if relaxing the above uniformality, it is highly desired that one may obtain a much sharper upper bound for $\fS(\mathbf{a};p)$ for many $\mathbf{a}\in\F_p^k.$ This is indeed the case when $\mathbf{a}\in\F_p^k$ runs over all possibilities as indicated by the equidistribution of Katz--Sarnak. Unfortunately, the range of $\mathbf{a}\in\F_p^k$ should be quite restrictive in our applications to Theorem \ref{thm:AB-bound-explicit}, and it seems quite difficult to recover a similar equidistribution when $\mathbf{a}$ runs over quite a thin family without friendly structures.

By closing this section, we would like to mention a dual but not quite related problem on behaviours of $\fS(\mathbf{a};p)$. Instead of varying $\mathbf{a}$ with a very large $p,$ we may alternatively consider the distribution of $\fS(\mathbf{a};p)$, with fixed $\mathbf{a}\in\Z^k,$ when $p$ runs over all large primes.
The case $k=4$ can be covered directly by the Sato--Tate conjecture for non-CM elliptic curves over $\Q$ as proven in a series of papers by L. Clozel, M. Harris and R. Taylor {\rm\cite{CHT08}}, R. Taylor {\rm\cite{Ta08}} and M. Harris, N. Shepherd-Barron and R. Taylor {\rm\cite{HST10}}, together with the CM case solved much earlier by E. Hecke {\rm\cite{He18}}. For $k\geqslant6$, one should refer to a suitable generalization of Sato--Tate conjecture in the framework of hyperelliptic curves; see {\rm\cite{KS09}} for instance.

\subsection*{Acknowledgments}
We are grateful to Daqing Wan for pointing out the reference
\cite{Wa97}. Y.-G. Chen is supported in part by the NSFC (No.
12171243), and P. Xi is supported in part by NSFC (No. 12025106,
No. 11971370) and The Young Talent Support Plan in Xi'an Jiaotong
University.

\smallskip


\smallskip

\begin{thebibliography}{abcdefg}

\bibitem[Bi68]{Bi68}
B.J. Birch,
How the number of points of an elliptic curve over a fixed prime field varies,
\emph{J. London Math. Soc.} \textbf{43} (1968) 57--60.



\bibitem[CY21]{CY21}
Y.-G. Chen and  X.-H. Yan, A conjecture of S\'ark\"ozy on
quadratic residues, \emph{J. Number Theory} \textbf{229} (2021)
100--124.


\bibitem[CHT08]{CHT08}
L. Clozel, M. Harris and R. Taylor,
Automorphy for some $\ell$-adic lifts of automorphic mod $\ell$ Galois representations, With Appendix A, summarizing unpublished work of Russ Mann, and Appendix B by Marie-France Vignéras,
\emph{Publ. Math. IH\'ES} \textbf{108} (2008) 1--181.

\bibitem[HST10]{HST10}
M. Harris, N. Shepherd-Barron and R. Taylor,
A family of Calabi--Yau varieties and potential automorphy,
\emph{Ann. of Math. (2)} \textbf{171} (2010) 779--813. 



\bibitem[He18]{He18}
E. Hecke, 
Eine neue Art von Zetafunktionen und ihre Beziehungen zur Verteilung der Primzahlen,
\emph{Math. Z.} \textbf{1} (1918) 357--376.


\bibitem[Ka92]{Ka92}
A.A. Karatsuba, The distribution of values of Dirichlet
characters on additive sequences, \emph{Soviet Math. Dokl.}
\textbf{44} (1992) 145--148.

\bibitem[KS99]{KS99}
N.M. Katz and P. Sarnak, 
Random matrices, Frobenius eigenvalues, and monodromy. 
AMS Colloquium Publications, No. 45, American
Mathematical Society, Providence, RI, 1999.

\bibitem[KS09]{KS09}
K.S. Kedlaya and A.V. Sutherland, 
Hyperelliptic curves, $L$-polynomials, and random matrices, 
Arithmetic, Geometry, Cryptography, and Coding Theory (AGCT 2007), \emph{Contemporary Math.} \textbf{487}, Amer. Math. Soc., 119-162, 2009.


\bibitem[Sa12]{Sa12}
A. S\'ark\"ozy, On additive decompositions of the set of quadratic
residues modulo $p$, \emph{Acta Arith.} \textbf{155} (2012)
41--51.

\bibitem[Shk14]{Shk14}
I.D. Shkredov, Sumsets in quadratic residues, \emph{Acta Arith.}
\textbf{164} (2014) 221--243.

\bibitem[Shp13]{Shp13}
I.E. Shparlinski, Additive decompositions of subgroups of finite
fields, \emph{SIAM J. Discrete Math.} \textbf{27} (2013)
1870--1879.


\bibitem[Ta08]{Ta08}
R. Taylor,
Automorphy for some $\ell$-adic lifts of automorphic mod $\ell$ Galois representations. II,
\emph{Publ. Math. IH\'ES} \textbf{108} (2008) 183--239. 




\bibitem[Wa97]{Wa97}
D. Wan, Generators and irreducible polynomials over finite fields,
\emph{Math. Comput.} \textbf{66} (1997) 1195--1212.

\bibitem[We48]{We48}
A. Weil, On some exponential sums, \emph{Proc. Nat. Acad. Sci.
U.S.A.} \textbf{34} (1948) 204--207.

\end{thebibliography}
\end{document}